\documentstyle[12pt]{article}
\if@twoside  m
\else
\fi
\newcommand{\ie}{{\em i.e.}}
\newcommand{\eg}{{\em e.g.}}
\newcommand{\rhs}{{\em rhs }}
\newcommand{\lhs}{{\em lhs }}
\newcommand{\OO}{{\cal O}}
\newcommand{\eps}{{\varepsilon}}
\newcommand{\R}{I\!\!R}

\newcommand{\QED}{\mbox{\rule[-1.5pt]{6pt}{10pt}}}
\newtheorem{claim}{Claim}[section]
\newtheorem{theorem}[claim]{Theorem}
\newtheorem{proposition}[claim]{Proposition}
\newtheorem{lemma}[claim]{Lemma}
\newtheorem{remark}[claim]{Remark}

\begin{document}
\vspace*{20mm}

\noindent
{\Large\bf Asymptotic estimates for bound states \\ in quantum
waveguides coupled laterally \\ through a narrow window}
\vspace{5mm}

\noindent
P.~Exner$^{1,2}$, S.A.~Vugalter$^{1,3}$ \\
{\small $^1$ Nuclear Physics Institute, Academy of Sciences,
25068 \v{R}e\v{z} near Prague, Czech Republic, \\ 
\phantom{$^1$}{\em exner@ujf.cas.cz} \\
$^2$ Doppler Institute, Czech Technical University, B\v rehov\'a 7,
11519 Prague, Czech Republic \\ 
$^3$ Radiophysical Research Institute, B. Pecherskaya 25/14, 603600
Nizhni Novgorod, Russia} 

\begin{quote}
{\bf Abstract.} Consider the Laplacian in a straight planar strip of
width $\,d\,$, with the Neumann boundary condition at a segment of
length $\,2a\,$ of one of the boundaries, and Dirichlet otherwise.
For small enough $\,a\,$ this operator has a single eigenvalue
$\,\epsilon(a)\,$; we show that there are positive $\,c_1,c_2\,$
such that $\,-c_1 a^4 \le \epsilon(a)- \left(\pi/ d\right)^2
\le -c_2 a^4\,$. An analogous conclusion holds for a pair of
Dirichlet strips, of generally different widths, with a window of
length $\,2a\,$ in the common boundary.
\end{quote}

\section{Introduction}

Recent progress in ``mesoscopic" physics brought not only new
physical effects but also some interesting spectral problems. 
One of them concerns the existence of bound states which appear if a
Dirichlet tube of a constant cross section is locally deformed, \eg,
bent or protruded, or coupled to another tube --- see 
\cite{BGRS,ES,DE,SRW} and references therein. 
In this paper we are concerned with another system of this type,
which consists of a pair of parallel Dirichlet strips coupled
laterally through a window in the common boundary; if they are of the
same width the problem simplifies to a treatment of a single strip
with the Dirichlet boundary condition changed to Neumann at a segment
of the boundary. 

Such a window--coupled system represents an idealized setup for some
existing quantum--wire devices \cite{HTW,Ku,LS,WG}. Its spectral and
scattering properties were discussed in recent papers \cite{BGRS,ESTV}.
It was shown there, in particular, that the discrete spectrum was
nonempty for any window width $\,2a>0\,$; if the latter is small
enough, there is just one simple eigenvalue $\,\epsilon(a)\,$ below
the bottom of the essential spectrum.

A question naturally arises about the behavior of the gap as
$\,a\to 0\,$. The result of \cite{BGRS} in combination with a
simple bracketing argument shows that it is bounded from below by
$\,c\,a^{4+\varepsilon}$; a numerical analysis performed in
\cite{ESTV} suggests that the asymptotic behavior is governed by the
fourth power of the window width. Our goal here is to prove two--sided
asymptotic estimates of this type.

\section{The results}

Consider a straight planar strip $\,\Sigma:= \R\times [-d_2,d_1]\,$.
Given $\,a>0\,$, we denote by $\,H(d_1,d_2;2a)\,$ the
Laplacian on $\,L^2(\Sigma)\,$ subject to the Dirichlet condition at
$\,y=-d_2,d_1\,$ as well as at the $\,|x|>a\,$ halfline segments of
the $\,x$--axis; this operator coincides with the Dirichlet Laplacian
at the strip with the cuts --- see \cite[Sec.XIII.15]{RS4}. Set
$\,d:=\max\{d_1,d_2\}\,$. If $\,d_1=d_2\,$, the operator decomposes
into an orthogonal sum with respect to the $\,y$--parity; the
nontrivial part is unitarily equivalent to the Laplacian on
$\,L^2(\Sigma_+)\,$, where $\,\Sigma_+:= \R\times [0,d]\,$, with the
Neumann condition at the segment $\,[-a,a]\,$ of the $\,x$--axis and
Dirichlet at the remaining part of the boundary; we denote it by
$\,H(d;2a)\,$.

Basic spectral properties of these operators are the following
\cite{BGRS,ESTV}: 

\begin{proposition}
$\sigma_{ess}(H(d_1,d_2;2a))= [(\pi/d)^2\!,\infty)$. The discrete
spectrum is non\-empty for any $\,a>0\,$ and consists of a finite
number of simple eigenvalues $\,(\pi/D)^2< \epsilon_1(a)<
\dots\epsilon_N(a) <(\pi/d)^2$, where $\,D:=d_1\!+d_2\,$; if
$\,a\,$ is small enough there is just one $\,\epsilon(a)\equiv
\epsilon_1(a)\,$.  
\end{proposition}

Our main result are asymptotic bounds for the narrow--window
case: 

\begin{theorem} \label{asymptotic bounds thm}
There are positive $\,c_1,\,c_2\,$ such that
   \begin{equation} \label{asymptotic bounds}
-c_1 a^4 \,\le\, \epsilon(a) -\left(\pi\over d\right)^2 \le -c_2
a^4 
   \end{equation}
holds for all sufficiently small $\,a\,$.
\end{theorem}

A proof of this theorem is the contents of the following sections.

\section{An upper bound}
\setcounter{equation}{0}

Let us begin with the operator $\,H\equiv H(d;2a)\,$. We denote by
$\,\chi_n(y)= \sqrt{2\over d}\, \sin\left(\pi ny\over d\right)\,$,
$\;n=1,2,\dots\,$, elements of the ``transverse" orthonormal basis
(the symbol should not be confused with the indicator function
$\,\chi_M\,$ of a set $\,M\,$), and set $\,\psi=F+G\,$, where 
   \begin{equation} \label{trial F}
f(x,y)\,:=\, f_1(x)\chi_1(y)\,, \qquad f_1(y)\,:=\, \alpha\,\max
\{\,1,\, e^{-\kappa||x|-a|}\}\,,
   \end{equation}
and
   \begin{equation} \label{trial G}
G(x,y)\,:=\,\eta\,\chi_{[-a,a]}(x)\, \cos\left(\pi x\over 2a\right)\,
R(y)
   \end{equation}
with
   \begin{equation} \label{R}
R(y)\,:=\, \left\lbrace\; \begin{array}{lcc} e^{-\pi y/2a} & \quad
\dots\quad & y\in\left\lbrack 0,{d\over2}\,\right\rbrack \\ \\
2\left(1-{y\over d}\,\right) e^{-\pi d/4a} & \quad
\dots\quad & y\in\left\lbrack {d\over2},d\,\right\rbrack \end{array}
\right. 
   \end{equation}
To make $\,\psi\,$ a trial function of $\,Q(H)\,$ we have to
ensure that it satisfies the Neumann boundary condition at the
window, \eg, replacing $\,G\,$ by $\,G^{\eps}\,$ such that
$$
G^{\eps}(x,d)=0\,,\qquad G^{\eps}_y(x,0)=\, -\, {\sqrt{2}\,
\alpha\pi \over d^{3/2}}\,,
$$
for $\,|x|\le a\,$, where the subscript is a shorthand for the partial
derivative, and
$$
\max\left\{\, \|G^{\eps}\!-G\|,\;  \|\nabla G^{\eps}\!-\nabla G\|,\;
\|G^{\eps}(\cdot,0)\!-G(\cdot,0)\|_{L^2(-a,a)}\: \right\} \,<\,\eps\,.
$$
We have to compute $\,L(\psi^{\eps}):= (H\psi^{\eps}, \psi^{\eps})-\,
\left(\pi\over d\right)^2\|\psi^{\eps}\|^2$. The operator part equals
$\,\|\psi^{\eps}_x\|^2\!+ \|\psi^{\eps}_y\|^2$, where the second term
is evaluated using $\,-\chi''_1= \left(\pi\over d\right)^2 \chi_1\,$,
a simple integration by parts, and the explicit value $\,\chi'_1(0)=
\pi\sqrt{2/d^3}\,$; together we get 
   \begin{equation} \label{trial L}
L(\psi^{\eps})\,=\, \|\psi^{\eps}_x\|^2+ \|G^{\eps}_y\|^2 -\,
\left(\pi\over d\right)^2 \|G^{\eps}\|^2- 2\alpha\, {\pi\over d}\,
\sqrt{2\over d}\, \int_{-a}^a G^{\eps}(x,0)\, dx\,.
   \end{equation}
It is sufficient to find the \rhs for $\,\psi\,$ since the difference
can be made arbitrarily small by a suitable choice of $\,\eps\,$.

Since $\,f_x,\, G_x\,$ have disjoint supports, we have
$\,\|\psi_x\|^2= \|F_x\|^2\!+ \|G_x\|^2$. Furthermore, the last term
equals $\,\eta^2\, {\pi^2\over 4a}\, \|R\|^2_{L^2(0,d)}\,$, and
$$
\|R\|^2_{L^2(0,d)}\,=\, {a\over\pi}\,+\, \left( {d\over 6}-
{a\over\pi} \right)\, e^{-\pi d/2a}\,<\, {a\over\pi}\,(1+\eps_1)
$$
for any $\,\eps_1>0\,$ and all $\,a\,$ small enough. In the same way,
$\,\int_{-a}^a G(x,0)\,dx=\, {4a\eta\over\pi}\,$. Finally, a bound to
$\,\|G_y\|^2$ follows from
$$
\|R'\|^2_{L^2(0,d)}\,=\, {\pi\over 4a}\,+\, \left( {2\over d}-
{\pi\over 4a} \right)\, e^{-\pi d/2a}\,<\, {\pi\over 4a}
$$
for $\,a<\pi d/8\,$, which means that $\,\|G_y\|^2< \eta^2\pi/4$.
Putting these estimates together, using $\,\|F_x\|^2=\alpha^2
\kappa\,$, and neglecting the negative term $\,-\left(\pi\over
d\right)^2 \|G\|^2\,$, we arrive at the inequality 
$$
L(\psi)\,<\, \alpha^2\kappa -\, {8\sqrt{2}\, \alpha a\over d^{3/2}}\,
\eta\,+\, {\pi\over 4}\, (2+\eps_1)\eta^2\,,
$$
where the sum of the last two terms at the \rhs is minimized by $\,-\,{2^7
\alpha^2 a^2\over \pi d^3(2+\eps_1)}\,$. It remains to estimate
$\,\|\psi\|^2$ from below. The tail contribution is
$\,\|\psi\|^2_{|x|\ge a}=\, \alpha^2/\kappa\,$, while the window part
expresses as 
$$
\|\psi\|^2_{|x|\le a} \,\le \, 2 \|F\|^2_{|x|\le a}+
2\|G\|^2_{|x|\le a} \,=\, 4a\alpha^2+ 2a\eta^2 \|R\|^2_{L^2(0,d)}
\,<\, Ca 
$$
for some $\,C>0\,$; this means that $\,\|\psi\|^2> \alpha^2
(1\!-\!\eps_2) \kappa^{-1}\,$ holds for any $\,\eps_2>0\,$ and all
$\,a\,$ small enough. Hence
$$
{L(\psi)\over\|\psi\|^2}\,<\, (1-\eps_2)^{-1} \left(\, \kappa^2-\,
{2^7 a^2\kappa\over \pi d^3(1+\eps_1)}\, \right)\;;
$$
taking the minimum over $\,\kappa\,$ we find
   \begin{equation} \label{trial bound}
{L(\psi)\over\|\psi\|^2}\,<\, (1-\eps_2)^{-1} {2^{12} a^4\over \pi^2
d^6(1+ \eps_1)^2}\,,
   \end{equation}
which completes the proof of the first inequality of (\ref{asymptotic
bounds}) in the symmetric case, $\,d_1=d_2\,$. 

Let us pass to the nonsymmetric case and suppose for definiteness
that $\,d=d_1>d_2\,$; the bottom of the essential spectrum is
then determined by the upper part of $\,\Sigma\,$. We choose there
the same trial function as above, \ie, $\,\psi=F+G\,$ is for $\,y\ge
0\,$ given by (\ref{trial F}) and (\ref{trial G}). Let further
$\,R_2\,$ be defined by (\ref{R}) with $\,d,y\,$ replaced by $\,d_2\,$
and $\,-y\,$, respectively. In the lower part of the strip,
$\,\R\times[-d_2,0]\,$, we put $\,\psi=G\,$ where $\,G\,$ is given by
(\ref{trial G}) with $\,R\,$ replaced by $\,R_2\,$. The trial
function should be smoothed in the window by requiring the
continuity, $\,G^{\eps}(x,0+)= G^{\eps}(x,0-)\,$, and 
$$
G^{\eps}_y(x,0+)\,=\, -\, {\sqrt{2}\,\alpha\pi \over d^{3/2}}\,,
\qquad G^{\eps}_y(x,0-)=0\,.
$$
Then it belongs to $\,Q(H)\,$; in the same way as above it is
sufficient to compute the functional for the nonsmoothed 
$\,\psi\,$. The value of $\,L(\psi)\,$ can be estimated by
$$
L(\psi)\,<\, \alpha^2\kappa -\, {8\sqrt{2}\, \alpha a\over d^{3/2}}\,
\eta\,+\, {\pi\over 2}\, (2+\eps_1)\eta^2\,,
$$
and since the contribution to $\,\|\psi\|^2$ from the window part is
again $\,\OO(a)\,$, the argument used in the symmetric case may be
repeated. \quad \QED

\section{Two lemmas}
\setcounter{equation}{0}

To prove the other part of Theorem~\ref{asymptotic bounds thm} we
need a few variational results for real--valued functions. Let us
first mention two elementary inequalities obtained directly by
solving the appropriate Euler equations. Let $\,\phi\in C^2(\R_+)\,$
with $\,\phi(0)=\alpha\,$; then  
   \begin{equation} \label{lemma 1}
\int_0^{\infty} \left(\phi'(t)^2\!+m^2\phi^2(t)\right)(t)\,dt \,\ge\,
m\alpha^2
   \end{equation}   
holds for a fixed $\,m\ge 0\,$. Similarly, if $\,\phi\in C^2[-b,b]$
with $\,\phi(\pm b)=0\,$, then 
   \begin{equation} \label{lemma 2}
\int_{-b}^b \phi'(t)^2 dt \,\ge\, \left(\pi\over 2b\right)^2
\int_{-b}^b \phi(t)^2 dt\,.
   \end{equation}

The following results are a little more involved:
  \begin{lemma} \label{lemma 3}
Let $\,\phi\in C^2[0,d]\,$ be a function with $\,\phi(d)=0\,$; then
there are positive $\,\eps_1,\eps_2\,$ such that $\,\left| \int_0^d
\phi(t)\chi_1(t)\,dt\right|< \eps_1\|\phi\|\,$ implies 
   \begin{equation} \label{derivative bound}
\int_{-d}^d \phi'(t)^2 dt\,>\, (1+\eps_2) \left(\pi\over d\right)^2
\|\phi\|^2\,.
   \end{equation}
   \end{lemma}
{\em Proof.} Let $\,\{g_n\}\,$ be the ``Dirichlet" trigonometric
basis in $\,L^2(-d,d)\,$; in particular, $\,2^{1/2} g_2\,$ is the
odd extension of $\,\chi_1\,$. We denote by $\,\Phi\,$ the even
extension of $\,\phi\,$ to $\,[-d,d]\,$; without loss of generality
we may suppose that it has a unit norm. Since $\,\Phi\,$ is even by
assumption, we write it as $\,\Phi= \gamma_1 g_1+h\,$ with
$\,h\in \{g_1,g_2\}^{\perp}$ which implies
   \begin{equation} \label{phi' norm}
\|\Phi'\|^2=\, \gamma_1^2 \left(\pi\over 2d\right)^2+ \|h'\|^2\,,
\qquad \|h'\|^2\ge\, \left(3\pi\over 2d\right)^2 \|h\|^2\,.
   \end{equation}
In a similar way we have $\,\xi=\xi_1g_1+\tilde\xi\,$ for the even
function $\,\xi:=|g_2|\,$. We find easily $\,\xi_1= 8/3\pi\,$, so
$\,\|\tilde\xi\|= \sqrt{1-\xi_1^2} \approx 0.529\,$.

If $\,|(\xi,\Phi)|<2^{3/2}\eps_1\,$, the identities $\,\xi_1\gamma_1=
(\xi,\Phi)- (\tilde\xi,h)\,$ and $\,\|h\|^2= 1-\gamma_1^2\,$ yield 
$$
{8\over 3\pi}\, |\gamma_1|\,\le\, \eps_1\sqrt{2}+\, \sqrt{1-\gamma_1^2}\,
\sqrt{1-\left(8\over 3\pi\right)^2}\,.
$$
This requires
$$
|\gamma_1|\,\le\, {8\sqrt{2}\over 3\pi}\,\eps_1\,+\,\sqrt{\left(1
-2\eps_1^2\right) \left(1-\left(8\over 3\pi\right)^2\right)} \;;
$$
hence choosing $\,\eps_1\,$ small enough we achieve $\,|\gamma_1|>\,
{1\over 2}\,$. Consequently, $\,\|\Phi'\|^2>\, {7\pi^2\over 4d^2}\,$
holds in view of (\ref{phi' norm}). \quad \QED

  \begin{lemma} \label{lemma 4}
Let $\,\phi\in C^2[0,d]\,$ with $\,\phi(0)=\beta\,$ and
$\,\phi(d)=0\,$. If $\,(\phi,\chi_1)=0\,$, then for every $\,m>0\,$
there is $\,c_0>0\,$ such that
   \begin{equation} \label{combined bound}
\int_0^d \phi'(t)^2 dt\,+\,\left(m\over a\right)^2 \int_0^a \phi(t)^2
dt\,-\, \left(\pi\over d\right)^2 \int_0^d \phi(t)^2 dt\,\ge\,
{c_0\beta^2 \over a}
   \end{equation}
holds for all $\,a\,$ small enough.
   \end{lemma} 
{\em Proof.} We denote the \lhs of (\ref{combined bound}) by
$\,M(\phi)\,$ and use subscripts to mark a contribution to a norm
from a particular interval. Furthermore, we introduce
   \begin{equation} \label{tilde phi}
\tilde\phi(x)\,:=\, \left\lbrace \begin{array}{ccc} \phi(a) & \quad
\dots \quad & 0\le x\le a \\ \phi(x) & \quad \dots \quad & a\le x\le
d \end{array} \right.
   \end{equation}
Notice that $\,\|\chi_1\|_{t\le a}\to 0\,$ as $\,a\to 0\,$, so there
is $\,a_0>0\,$ such that $\,2\|\chi_1\|_{t\le a}<\eps_1\,$ holds for
$\,a<a_0\,$, where $\eps_1\,$ is the positive number from the
previous lemma; we shall restrict ourselves in the following to
   \begin{equation} \label{a restriction}
a\,<\, \min\left\lbrace\, a_0,\, {md\over \pi\sqrt 3}
\,\right\rbrace\,. 
   \end{equation}
Suppose first that 
   \begin{equation} \label{nonuniform approach}
\|\phi\|_{t\le a}^2 >\, \|\tilde\phi\|^2 
   \end{equation}
holds for some $\,\phi\,$. Since the \rhs is not smaller than
$\,\|\phi\|_{t\ge a}^2$, we have in view of (\ref{a restriction})
$$
\max\left\lbrace\, \left(\pi\over d\right)^2\|\phi\|^2_{t\le a},\,
\left(\pi\over d\right)^2\|\phi\|^2_{t\ge a} \right\rbrace \,<\,
{m^2\over 3a^2}\,\|\phi\|^2_{t\le a}\,,
$$
so neglecting the non--negative term $\,\|\phi'\|^2_{t\ge a}$, we
arrive at the estimate
   \begin{equation} \label{nonuniform estimate}
M(\phi)\,>\, \|\phi'\|^2_{t\le a}+\, {m^2\over 3a^2}\,
\|\phi\|^2_{t\le a}\,\ge\, {\beta^2 m\over a\sqrt 3}\, \tanh {m
\over\sqrt 3}\;;
   \end{equation}
the last inequality follows from the fact that the estimating
functional is for a fixed $\,\alpha:= \phi(a)\,$ minimized by
$\,\phi(t)= {\beta\over\sinh\mu}\, \left(\alpha\sinh \mu t+ \sinh\mu
(a\!-\!t)\right)\,$, where $\,\mu:= {m\over\sqrt 3}\,$, and taking
the minimum over $\,\alpha\,$. Consider on the contrary those
$\,\phi\,$ for which (\ref{nonuniform approach}) is violated, \ie,
   \begin{equation} \label{uniform approach}
\|\phi\|_{t\le a}^2 \le\, \|\tilde\phi\|^2\,.
   \end{equation}
Since $\,\phi\,$ is by assumption orthogonal to $\,\chi_1\,$, we have
$$
\int_0^d \tilde\phi(t)\chi_1(t)\, dt \,=\, \int_0^a
\tilde\phi(t)\chi_1(t)\, dt \,-\, \int_0^a
\phi(t)\chi_1(t)\, dt\,,
$$
so
$$
\left|\int_0^d \tilde\phi(t)\chi_1(t)\, dt\,\right| \,\le\,
\left(\|\tilde\phi\|_{t\le a}+ \|\phi\|_{t\le a} \right)
\|\chi_1\|_{t\le a} \,\le\, 2\|\chi_1\|_{t\le a}
\|\tilde\phi\|\,<\,\eps_1\|\tilde\phi\|\;;
$$
applying Lemma~\ref{lemma 3} to $\,\tilde\phi\,$ we obtain 
$$
\int_a^d \phi'(t)^2 dt\,-\, \left(\pi\over d\right)^2 \int_a^d
\phi(t)^2 dt \,\ge\, \int_0^d \tilde\phi'(t)^2 dt\,-\, \left(\pi\over
d\right)^2 \int_a^d \tilde\phi(t)^2 dt\,>\, \eps_2 \left(\pi\over
d\right)^2 \|\tilde\phi\|^2
$$
for some $\,\eps_2>0\,$. Neglecting this non--negative term, we get
   \begin{equation} \label{uniform estimate}
M(\phi)\,>\, \|\phi'\|^2_{t\le a}+\, {m^2\over a^2}\,
\|\phi\|^2_{t\le a}\,-\,\left(\pi\over d\right)^2 \|\phi\|_{t\le
a}^2\;;
   \end{equation}
then it is sufficient to employ (\ref{a restriction}) and estimate
the \rhs in the same way as in (\ref{nonuniform estimate}) to arrive
at the desired conclusion.  
\quad \QED

\section{A lower bound}
\setcounter{equation}{0}

This is the most difficult part of the proof. However, we may
restrict ourselves to the symmetric case only because inserting an
additional Neumann boundary into the window we get a lower bound to
$\,H(d_1,d_2;2a)\,$; hence it is sufficient to treat the spectrum
of $\,H\equiv H(d;2a)\,$.  

We have to estimate $\,L(\psi):= (H\psi, \psi)-\, \left(\pi\over
d\right)^2\|\psi\|^2$ over $\,\|\psi\|^2$ for all $\,\psi\,$ of a
core of $\,H$, say, all $\,C^2$--smooth $\,\psi\in L^2(\Sigma_+)\,$
satisfying the boundary conditions. We shall employ for such
functions and $\,|x|\ge a\,$ the following uniformly convergent
Fourier expansion 
   \begin{equation} \label{Fourier}
\psi(x,y)\,=\, \sum_{n=1}^{\infty} c_n(x)\chi_n(y)\,,
   \end{equation}
where the coefficients $\,c_n(x)=(\psi(x,\cdot),\chi_n)\,$ are again
smooth. Some natural restrictions may be adopted:
   \begin{description}
   \vspace{-.8ex}
\item{\em (i)} only real--valued $\,\psi\,$ should be taken into
account: since $\,H\,$ commutes with complex conjugation we may
consider $\,L^2(\Sigma_+)\,$ as a real Hilbert space in which
$\,H\,$ is ``doubled",
   \vspace{-1.8ex}
\item{\em (ii)} we may consider $\,x$--even $\,\psi\,$ only,
because due to mirror symmetry we have $\,H=H_{even}\oplus H_{odd}$,
where the two parts are unitarily equivalent to the halfstrip
Laplacian with the Neumann and Dirichlet condition, respectively, at
the transverse cut. Hence $\,H_{even}\le H_{odd}$, and it is sufficient
to estimate the even part only. 
   \vspace{-.8ex}
   \end{description}
Next we have to introduce some more notation. We put $\,c_1=f_1\!
+\hat f_1\,$, where 
   \begin{equation} \label{f_1}
\hat f_1\,:=\, \left\lbrace \begin{array}{lll} 0 & \quad \dots \quad
& |x|\ge 2a \\ c_1(x)\!-\!\alpha & \quad \dots \quad & |x|\ge 2a
\end{array} \right.
   \end{equation}
with $\,\alpha:=c_1(2a)\,$, so $\,\hat f_1(\pm 2a)=0\,$, and
furthermore 
   \begin{equation}
\psi(x,y) \,=\, F(x,y)+G(x,y)\,, \qquad F(x,y)\,:=\,
f_1(x)\chi_1(y) \label{psi}\,.
   \end{equation}
Using this decomposition we derive in the same way as in Section~3 an
expression for the functional to be estimated,
   \begin{equation} \label{L}
L(\psi)\,=\, \|\psi_x\|^2+ \|G_y\|^2 -\,
\left(\pi\over d\right)^2 \|G\|^2- 2\alpha\, {\pi\over d}\,
\sqrt{2\over d}\, \int_{-a}^a G(x,0)\, dx\,.
   \end{equation}
Let us begin with the contribution from the outer region. We split a
half of the first term and consider the following expression:
   \begin{eqnarray*}
\lefteqn{{1\over 2}\,\|\psi_x\|^2_{|x|\ge a}+ \|G_y\|^2_{|x|\ge a}
-\,\left(\pi\over d\right)^2 \|G\|^2_{|x|\ge a}} \\ \\ && =\,
\sum_{n=1}^{\infty} \int_a^{\infty} c_n'(x)^2\, dx
\,+\,2\sum_{n=1}^{\infty} \left(\pi\over d\right)^2
\left(n^2\!-1\right) \int_a^{\infty} c_n(x)^2\, dx \,. 
   \end{eqnarray*}
Since 
$$
\int_a^{\infty} \left\lbrack\, {1\over 2}\, c_n'(x)^2\, dx +
\left({\pi\over d}\,\sqrt{n^2\!-1}\,\right)^2 c_n(x)^2\,\right\rbrack\,
dx\,>\, {\pi n\over 2d}\, c_n(a)^2
$$
by (\ref{lemma 1}) for $\,n\ge 2\,$, and the non--negative term
$\,\int_a^{\infty} c_1'(x)^2\, dx\,$ may be neglected, we have 
$$
{1\over 2}\,\|\psi_x\|^2_{|x|\ge a}+ \|G_y\|^2_{|x|\ge a}
-\,\left(\pi\over d\right)^2 \|G\|^2_{|x|\ge a}\,>\, {\pi\over d}\,
\sum_{n=2}^{\infty} nc_n(a)^2\,,
$$
and therefore
   \begin{eqnarray} \label{L1}
L(\psi) &\!>\!& {1\over 2}\,\|\psi_x\|^2_{|x|\ge a}+ \|G_x\|^2_{|x|\le
a} + \|G_y\|^2_{|x|\le a} -\,\left(\pi\over d\right)^2
\|G\|^2_{|x|\le a} \nonumber \\ \\ &\!+\!&{\pi\over d}\,
\sum_{n=2}^{\infty} nc_n(a)^2 -\, 2\alpha\, {\pi\over d}\,
\sqrt{2\over d}\, \int_{-a}^a G(x,0)\, dx\,, \nonumber
   \end{eqnarray}
where we have also employed the fact that $\,\psi_x=G_x\,$ for
$\,|x|\le 2a\,$. 

Our next goal is to estimate the sum of the first two terms in
(\ref{L1}) by estimating $\,\|G_x\|^2_{|x|\le 2a}\,$. If a function
$\,\tilde G:\, \Sigma_+\to C^2(\Sigma_+)\,$ vanishes at $\,x=\pm 2a\,$,
the inequality (\ref{lemma 2}) implies 
   \begin{equation} \label{zero end bound}
\|\tilde G_x\|^2_{|x|\le 2a}\,\ge\, \left( \pi\over 4a\right)^2
\|\tilde G\|^2_{|x|\le 2a}\,.
   \end{equation}
Unfortunately, $\,G\,$ does not satisfy this requirement, which
forces us to split it into several components and to estimate them
separately. First we single out the projection of $\,G\,$ onto the
first transverse mode,
   \begin{equation} \label{G}
G(x,y) \,=\, G_1(x,y)+G_2(x,y)\,, \qquad G_1(x,y)\,=\,\hat
f_1(x)\chi_1(y)\,.
   \end{equation}
Since $\,\hat f_1(\pm 2a)=0\,$ by definition, (\ref{zero end bound})
may be applied to $\,G_1\,$. In combination with the inequality
$$
{1\over 2}\,\|\psi_x\|^2_{a\le |x|\le 2a}+\,\|G_x\|^2_{|x|\le a}
\,\ge\, {1\over 2}\,\|G_x\|^2_{|x|\le 2a} \,=\, {1\over
2}\,\|G_{1,x}\|^2_{|x|\le 2a}+\, {1\over 2}\,\|G_{2,x}\|^2_{|x|\le
2a}\,, 
$$
where the subscript in the last norm may be dropped because
$\,G_1(x,y)=0\,$ for $\,|x|\ge 2a\,$, we get
   \begin{eqnarray} \label{L2}
L(\psi) &\!>\!& {1\over 2}\,\|\psi_x\|^2_{|x|\ge 2a}+\,{1\over 2}\,
\|G_{2,x}\|^2_{|x|\le 2a} + \|G_y\|^2_{|x|\le a} -\,\left(\pi\over
d\right)^2 \|G\|^2_{|x|\le a} \phantom{MMM} \nonumber \\ \\
&\!+\!&{\pi\over d}\, \sum_{n=2}^{\infty} nc_n(a)^2 -\, 2\alpha\,
{\pi\over d}\, \sqrt{2\over d}\, \int_{-a}^a G(x,0)\, dx \,+\,
{\pi^2\over 32 a^2}\, \|G_1\|^2\,. \nonumber
   \end{eqnarray}
The function $\,G_2\,$ has again to be splitted; we rewrite it as 
   \begin{equation} \label{G_2}
G_2(x,y) \,=\, \hat G(x,y)+\Gamma(x,y)\,, \qquad \Gamma(x,y)\,:=\,
\sum_{n=2}^{\infty} c_n(2a)\chi_n(y)\,,
   \end{equation}
for $\,|x|\le 2a\,$; the second part is independent of $\,x\,$ while
the first one vanishes at the border, $\,|x|=2a\,$, so $\,G_{2,x}=
\hat G_x\,$ may be estimated by means of (\ref{zero end bound}) and
the Schwarz inequality as 
   \begin{equation} \label{G_2 bound}
\|G_{2,x}\|^2_{|x|\le 2a}\ge\, \left(\pi\over 4a\right)^2 \|\hat
G\|^2_{|x|\le 2a} \,\ge\, {1\over 2}\left(\pi\over 4a\right)^2 \|
G_2\|^2_{|x|\le 2a}\!- \left(\pi\over 4a\right)^2 \|\Gamma
\|^2_{|x|\le 2a}.   
   \end{equation}
However, it is difficult to find a suitable bound for the last
negative term. Instead of attempting it we restrict ourselves to the
vicinity of the window: we introduce $\,\Omega_a:= [-2a,2a]\times
[0,a]\,$ and replace (\ref{G_2 bound}) by
   \begin{equation} \label{modified G_2 bound}
\|G_{2,x}\|^2_{|x|\le 2a}\ge\, \left(\pi\over 4a\right)^2 \|\hat
G\|^2_{|x|\le 2a} \,\ge\, \left(\pi\over 4a\right)^2 \|\hat
G\|^2_{\Omega_a} \,\ge\, {1\over 2}\left(\pi\over 4a\right)^2 \|
G_2\|^2_{\Omega_a}\!- \left(\pi\over 4a\right)^2 \|\Gamma
\|^2_{\Omega_a}. 
   \end{equation}
We have still to find an upper bound to $\,\|\Gamma\|^2
_{\Omega_a}\,$. To this end we notice that, in addition to 
{\em (i), (ii)}, we may restrict our attention to those $\,\psi\,$
for which
   \begin{equation} \label{subexponential decay}
|c_n(x)|\,\le\, |c_n^{ex}(x)|\,, \qquad c_n^{ex}(x)\,:=\, c_n(a)\,  
e^{-(\pi/d)\sqrt{n^2-1} (x-a)}\,,
   \end{equation}
holds for $\,|x|\ge a\,$ and $\,n\ge 2\,$. Indeed, let us split
$\,\psi\,$ for $\,|x|\ge a\,$ into a contribution from the $\,n$--th
transverse mode, $\,n\ge 2\,$, and the orthogonal complement
introducing 
$$
\tilde\psi(x,y)\,:=\, \left\lbrace \begin{array}{lll} \psi(x,y) &
\quad \dots \quad & |x|<a \\ \psi(x,y)-c_n(x)\chi_n(y) & \quad \dots
\quad & |x|\ge a \end{array} \right.
$$
The basic expression to be estimated can be then written as
$$
{L(\psi)\over \|\psi\|^2} \,=\,{\left(H\tilde\psi,
\tilde\psi\right)-\,\left(\pi\over d\right)^2 \|\tilde\psi\|^2+\,
2\int_a^{\infty} \left\lbrack\, c_n'(x)^2\, dx + \left({\pi\over
d}\,\sqrt{n^2\!-1}\,\right)^2 c_n(x)^2\,\right\rbrack\, dx  \over
\|\tilde\psi\|^2 +2 \int_a^{\infty} c_n(x)^2 dx} \,.
$$
Without loss of generality we may assume only those $\,\psi\,$ for
which the numerator is negative. The part of its last term
corresponding to the interval $\,[a,\infty)\,$ is by (\ref{lemma
1}) minimized by the exponential function $\,|c_n^{ex}|$ of
(\ref{subexponential decay}); hence replacing $\,c_n(x)^2\,$ by
$\,\min\{c_n(x)^2, c_n^{ex}(x)^2\}\,$ we can only get a larger
negative number. At the same time, the positive denominator can be
only diminished, which justifies the claim made above.

We are interested in the norm of $\,\Gamma\,$ restricted to
$\,\Omega_a\,$, hence we cannot use the Parseval relation because in
general the restrictions of $\,\chi_n\,$ to $\,[0,a]\,$ are not
orthogonal. We divide therefore the series into two pieces referring
to small and large values, respectively, relative to $\,a^{-1}$. In
the first part we employ the smallness of the $\,\chi_n\,$ norm
restricted to $\,[0,a]\,$, while the other part will be estimated by
means of the subexponential decay (\ref{subexponential decay}). In
this way, we may write
   \begin{eqnarray*}
\|\Gamma\|^2_{\Omega_a} &\!=\!& \int_{-2a}^{2a} dx \int_0^a dy\,
\left( \sum_{n=2}^{\infty} c_n(2a)\chi_n(y) \right)^2 \\ \\
&\!\le\!& 8a \int_0^a \left( \sum_{n=2}^{[a^{-1}]+1} c_n(2a)\chi_n(y)
\right)^2\! dy \,+\, 8a \int_0^a \left( \sum^{\infty}_{2\le
n=[a^{-1}]+2} c_n(2a)\chi_n(y) \right)^2 \!dy \\ \\
&\!\le\!& 8a \left( \sum_{n=2}^{[a^{-1}]+1} n^{-1} c_n(a)^2 \int_0^a
\chi_n(y)^2 dy \right) \left( \sum_{n=2}^{[a^{-1}]+1} n \right) \\ \\
&\!+\!& 8a \left( \sum^{\infty}_{2\le n=[a^{-1}]+2} n c_n(a)^2 \int_0^a
\chi_n(y)^2 dy \right) \left( \sum^{\infty}_{2\le n=[a^{-1}]+2} n^{-1}
e^{-(2\pi a/d)\sqrt{n^2-1}} \right)\,,
   \end{eqnarray*}
where $\,[\cdot]\,$ denotes the entire part; in the first term on the
\rhs we have used the rough bound $\,c_n(2a)<c_n(a)\,$. The
transverse--mode integral equals 
$$
\int_0^a \chi_n(y)^2 dy\,=\, {a\over d}\, \left\lbrack\,
1-\left(d\over 2\pi na\right)\sin\left(2\pi na\over d\right)
\,\right\rbrack \,\le\, {a\over d}\,\min\left\lbrace\, {1\over 6}\,
\left(2\pi na\over d\right)^2\!,\: 2\:\right\rbrace\;;
$$
being small for small $\,a\,$ as indicated. We arrive thus at the
inequality 
   \begin{equation} \label{Gamma bound}
\|\Gamma\|^2_{\Omega_a}\,\le\, {16a^2\over d}\, \left\lbrace\, 
{2\pi^2\over 3d^2}\,+\, \sum^{\infty}_{n=[a^{-1}]+1} n^{-1} e^{-2\pi
an/d} \,\right\rbrace\: \sum_{n=2}^{\infty} n c_n(a)^2\,,
   \end{equation}
where we have employed the estimates $\,\sum_{n=2}^{[a^{-1}]+1} n\,
\le 2a^{-1}\,$ and $\,\sqrt{n^2\!-1}< n\!-\!1\,$. The sum in the
curly bracket at the \rhs of (\ref{Gamma bound}) has a uniform upper
bound with respect  to $\,a\,$ being a Darboux sum of the integral
$$
\int_1^{\infty} {e^{-2\pi x/d}\over x}\: dx\,=\, -\,{\rm Ei}\left(
-\,{2\pi\over d}\right)\,,
$$
where $\,{\rm Ei}\,$ is the exponential integral function. Hence
there is a positive $\,C\,$ such that 
   \begin{equation} \label{Gamma bound 2}
\|\Gamma\|^2_{\Omega_a}\,\le\, Ca^2\, \sum_{n=2}^{\infty} n
c_n(a)^2\,. 
   \end{equation}
By (\ref{modified G_2 bound}) we have
$$
{1\over 2}\,\|G_{2,x}\|^2_{|x|\le 2a}+\, {\pi\over
d}\,\sum_{n=2}^{\infty} n c_n(a)^2 \ge\, {\delta\over
4}\left(\pi\over 4a\right)^2 \| G_2\|^2_{\Omega_a}\!- {\delta\over
2}\: \|\Gamma \|^2_{\Omega_a}+\, {\pi\over d}\,\sum_{n=2}^{\infty} n
c_n(a)^2 
$$
for any $\,\delta\in(0,1]\,$; the estimate (\ref{Gamma bound 2}) shows
that choosing $\,\delta\,$ small enough one can achieve that the sum
of the last two terms is non--negative. Putting $\,m:=\,
{\pi\over 8}\,\sqrt{\delta}\,$, we arrive at the bound
   \begin{eqnarray} \label{L3}
L(\psi) &\!>\!& {1\over 2}\,\|\psi_x\|^2_{|x|\ge 2a}+
\|G_y\|^2_{|x|\le a} -\,\left(\pi\over d\right)^2 \|G\|^2_{|x|\le a}
+\, {\pi^2\over 32 a^2}\, \|G_1\|^2 \phantom{MM} \nonumber \\ \\
&\!-\!& 2\alpha\,{\pi\over d}\, \sqrt{2\over d}\, \int_{-a}^a
G(x,0)\, dx \,+\, {m^2\over a^2}\, \|G_2\|^2_{\Omega_a}\,. \nonumber
   \end{eqnarray}

In the next step we express the term containing $\,G_y\,$ using the
decomposition (\ref{G}), a simple integration by parts, the relation
$\,G_2(x,0)= G(x,0)\,$, and the fact that $\,G_2(x,\cdot)\,$ is
orthogonal to $\,\chi''_1= -(\pi/d)^2\chi_1\,$. This yields
$$
\|G_y\|^2_{|x|\le a}=\, \|G_{1,y}\|^2_{|x|\le a}+\,
\|G_{2,y}\|^2_{|x|\le a} -\,2\,{\pi\over d}\, \sqrt{2\over d}\,
\int_{-a}^a \hat f_1(x)G(x,0)\, dx\,,
$$
where the last term does not exceed $\,(\pi/d^2)
\left(2\|G_1\|^2_{|x|\le a}+ d\|G(\cdot,0)\|^2_{|x|\le a} \right)\,$
by the Schwarz inequality. We substitute to (\ref{L3}) from here
keeping in mind the identity $\,\|G\|^2_{|x|\le a}= \|G_1\|^2_{|x|\le
a}\!+ \|G_2\|^2_{|x|\le a}\,$, and neglect the term
$\,\|G_{1,y}\|^2_{|x|\le a}\,$ as well as 
$$
{\pi^2\over 32 a^2}\, \|G_1\|^2-\, {\pi(\pi+2)\over d^2}\,
\|G_{1,y}\|^2_{|x|\le a} 
$$
which is positive for $\,a\,$ small enough, obtaining
   \begin{eqnarray} \label{L4}
L(\psi) &\!>\!& {1\over 2}\,\|\psi_x\|^2_{|x|\ge 2a}+\,
\|G_{2,y}\|^2_{|x|\le a} +\,{m^2\over a^2}\, \|G_2\|^2_{\Omega_a}
-\,\left(\pi\over d\right)^2 \|G_2\|^2_{|x|\le a}
\phantom{MMM}\nonumber \\ \\ &\!-\!& {\pi\over
d}\,\|G_2(\cdot,0)\|^2_{|x|\le a} -\,2\alpha\,{\pi\over d}\,
\sqrt{2\over d}\, \int_{-a}^a  G_2(x,0)\, dx \,. \nonumber
   \end{eqnarray}
Furthermore, the function $\,G_2(x,\cdot)\,$ satisfies for a fixed
$\,x\in [-a,a]\,$ the assumptions of Lemma~\ref{lemma 4}, so the sum
of the second, third, and fourth term on the \rhs is below bounded by
$\,(c_0/ a)\, \|G_2(\cdot,0)\|^2_{|x|\le a}\,$.  Since $\,(c_0/2a)
-(\pi/d)>0\,$ for small $\,a\,$, we have 
$$
L(\psi) \,>\, {1\over 2}\,\|\psi_x\|^2_{|x|\ge 2a}+\, {c_0\over 2a}\,
\|G_2(\cdot,0)\|^2_{|x|\le a} -\,2\alpha\,{\pi\over d}\,
\sqrt{2\over d}\,\|G_2(\cdot,0)\|_{|x|\le a} \sqrt{2a}\,,
$$
where we have employed the Schwarz inequality again; taking the
minimum over $\,\|G_2(\cdot,0)\|_{|x|\le a}$, we arrive finally at
the estimate 
   \begin{equation} \label{L5}
L(\psi) \,>\, {1\over 2}\,\|\psi_x\|^2_{|x|\ge 2a}
-\,{8\pi^2\alpha^2\over c_0d^3}\:.
   \end{equation}

The rest of the argument is simple. We have $\,\|\psi\|^2\ge
\|\psi\|^2_{|x|\ge 2a} \ge 2\int_{2a}^{\infty} c_1(x)^2 dx\,$, and a
similar bound is valid for the first term on the \rhs of (\ref{L5}),
so 
$$
{L(\psi)\over \|\psi\|^2} \,>\, {\int_{2a}^{\infty} c'_1(x)^2 dx
\,-\,{8\pi^2\alpha^2\over c_0d^3} \over 2\int_{2a}^{\infty} c_1(x)^2
dx}\;.
$$
The extremal of this functional over functions with a fixed value at
$\,x=2a\,$ is $\,c_1(2a)\, e^{-\kappa(x-2a)}$ which yields the value
$\,(\kappa^2/2)-(8\pi^2/c_0d^3)\kappa a^2\,$. It is now sufficient to
take the minimum over $\,\kappa\,$ to get the inequality
   \begin{equation} \label{lower bound}
{L(\psi)\over \|\psi\|^2} \,>\, -\,2\, \left( 4\pi^2\over c_0d^3
\right)^2 a^4\,,
   \end{equation}
which completes the proof of Theorem~\ref{asymptotic bounds thm}.
\quad \QED

\begin{remark}
{\rm Some of the estimates we used are certainly crude, however,
pushing them to an optimum is not sufficient to squeeze the bound
(\ref{asymptotic bounds}) enough to get the actual asymptotic
behavior. This question remains open; the same is true for the
window--width threshold behavior of higher eigenvalues \cite{ESTV}. }
\end{remark}
\vspace{3mm}

\noindent
{\em Acknowledgments.} S.V. is grateful for the hospitality
extended to him in the Nuclear Physics Institute, AS, where this work
was done. The research has been partially supported by the Grants AS
No. 148409, GACR No. 202--0218, and GosComVUZ 95--0--2.1--50.

\end{document}